\documentclass[10pt]{amsart}
\usepackage{amssymb,latexsym,amsmath,graphics}
\theoremstyle{plain}
\newtheorem{theorem}{Theorem}[section]
\newtheorem{corollary}{Corollary}[section]
\newtheorem{lemma}{Lemma}[section]
\newtheorem{definition}{Definition}[section]
\newtheorem*{example}{Example}

\begin{document}
\title{A Bijection Between Certain Non-Crossing Partitions and Sequences}
\author{Rekha Natarajan}
\maketitle
\begin{abstract}
We present a bijection between
non-crossing partitions of the set $[2n+1]$ into $n+1$ blocks such
that no block contains two consecutive integers, and the set of sequences
$\{s_{i}\}_{1}^{n}$ such that $1 \leq s_{i} \leq i$, and if $s_{i}=j$,
then $s_{i-r} \leq j-r$ for $1 \leq r \leq j-1$.
\end{abstract}
\section{Introduction}
The Catalan numbers, $C_{n}=\frac{1}{n+1}\binom{2n}{n}$, arise in
a wide selection of combinatorial problems and applications. In
his book \emph{Enumerative Combinatorics}, Volume II, R.P. Stanley
gives a list of 66 sets whose elements are counted by the Catalan
numbers \cite{Stan}. Later on, he has also provided an additional
26 sets in his Catalan addendum \cite{addendum}. Out of this vast
collection, we focus our attention on a specific type of
non-crossing partitions and sequences.

Non-crossing partitions, NCP's, of the set $[n]=\{1,\dots,n\}$,
where $n \in \mathbb{N}$, are partitions $\pi = \{B_{1}, \dots,
B_{k}\} \in \Pi_{n}$, such that if $a< b < c < d$, with $a,c \in
B_{i}$ and $b,d \in B_{j}$, then $i=j$ (here, $\Pi_{n}$ denotes
the set of all partitions of $[n]$.) Such partitions will always
be given in lexicographic order.  For instance when $n=5$,
$\pi=\{1,5\}-\{2,4\}-\{3\}$ is a non-crossing partition, whereas
$\pi=\{1,3\}-\{2,4\}-\{5\}$ is a crossing partition. We choose to
examine a specific class of non-crossing partitions, namely the
set of non-crossing partitions of the integer $[2n+1]$ into $n+1$
blocks such that no block contains two consecutive integers,
hereby referred to as $\emph{\textquoteleft special'}$ partitions
(\cite{Stan}, p. 226, Exercise 6.19(tt)).  In \cite{Mullin},
Mullin and Stanton prove that the set of special partitions of the
set $[2n+1]$ is in one-to-one correspondence with the set of plane
trees on $n+1$ vertices, which in turn has cardinality $C_{n}$
(see \cite{Stan}, p. 176). For more information on these trees,
see \cite{Klaner}. Other known bijections involving special
partitions can also be found in Roselle \cite{Roselle}.

In addition to special partitions, the Catalan numbers also count
sequences $s_{1}, s_{2}, \dots, s_{n}$ of integers such that $1
\leq s_{i} \leq i,$ and if $s_{i}=j$, then $s_{i-r} \leq j-r$ for
$1 \leq r \leq j-1$ (see \cite{Stan}, p. 223, Exercise 6.19(z)).
Stanley gives bijections between this set of sequences and the set
of $312$-avoiding permutations of $[n],$ as well as binary trees
on $n$ vertices (\cite{Stan}, p. 259, Exercise 6.19(z)).

In this paper, we present
a bijection between special partitions and the sequences described above.
Up until now, no bijections have been discovered between
these two sets.  We hope that this connection between
a certain class of non-crossing partitions and sequences will
shed a new light on the vast area of research encompassing
permutations with some restricted patterns.

Section 2 focuses on background information relating to special
partitions,
while section 3 discusses the main result, which is the bijection itself.

\section{Definitions and Results}

We first present a few definitions and preliminary results concerning special
partitions.  Some of the following results may already be known,
but we have not found any information regarding these in the literature.
In any case, we include them here for completeness.

\begin{definition}
Let $\mathcal P_{[2n+1]}$ be the set of non-crossing partitions of
$[2n+1]$ into $n+1$ parts such that no block contains two
consecutive integers. Any partition in the set $\mathcal
P_{[2n+1]}$ will be referred to as a $\textquoteleft$special'
partition.
\end{definition}
An element $\pi \in \mathcal P_{[13]}$, for example, is:
\begin{center}
    $\pi = \{1,13\} -
\{2,4,6,12\} - \{3\} - \{5\} - \{7,11\} - \{8,10\} -  \{9\}$
\end{center}
We will represent $\pi$ in the following manner (also known as a
Puttenham diagram \cite{Becker}).
\begin{figure}[ht]
\centering
    \scalebox{0.5}{\includegraphics{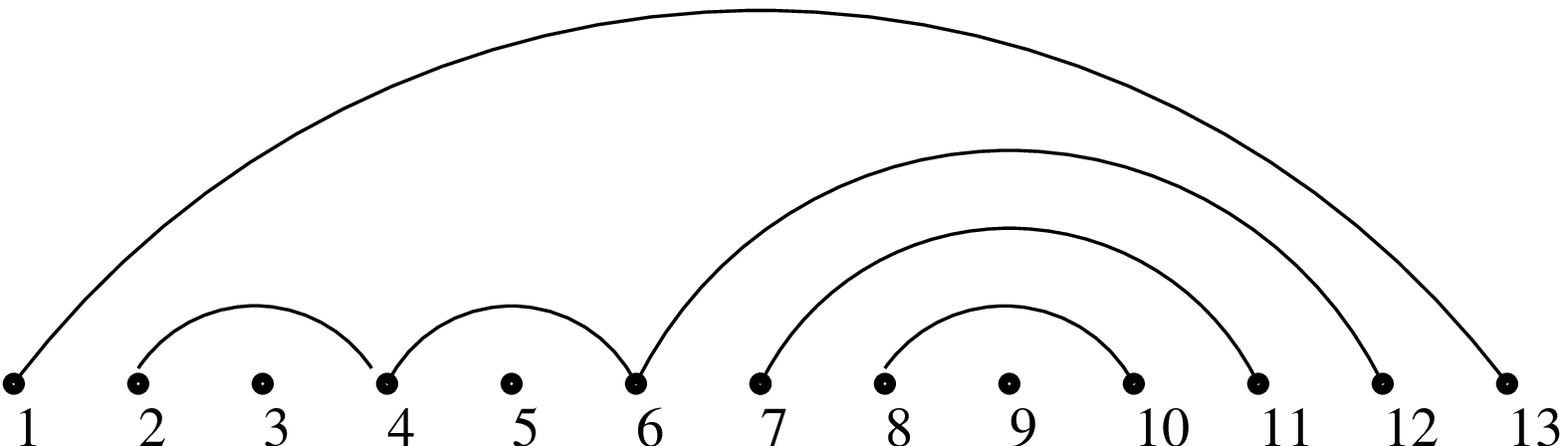}}
    \caption{$\pi = \{1,13\} - \{2,4,6,12\} - \{3\} - \{5\} - \{7,11\}
 - \{8,10\}
 -  \{9\}$}
\end{figure}
\begin{definition}
Given a non-crossing partition $\pi=\{B_{1}, B_{2}, \dots, B_{k}\}
\in \Pi_{n}$, we partition $\pi$ into (disjoint) \textbf{pieces}
in the following manner:
\begin{enumerate}
    \item If $B_{1}={1}$, then $B_{1}$ is said to be a piece of $\pi$.

    \item If $B_{1}=\{1,\dots,b_{1_{w}}\}$, let
$B_{2}, \dots, B_{1+s}$ be the set of blocks in $\pi$ containing
the integers $2,\dots,b_{1_{w}}-1$. Then the set of blocks $B_{1},
B_{2}, \dots, B_{1+s}$ is said to be a piece of $\pi$.\\
To find the next piece, consider
$B_{1+s+1}=B_{t}=\{b_{1_{w}}+1=b_{t_{1}}, \dots, b_{t_{m}}\}$.
Let $B_{t+1},\dots,B_{t+i}$ be the set of blocks containing
$b_{t_{1}}+1, \dots, b_{t_{m}}-1$. Then the
set of blocks $B_{t}, B_{t+1}, \dots, B_{t+i}$ is also said to be
a piece of $\pi$.  Note that this piece could also
be a singleton. Continue this process
until all elements of $[n]$ have been considered.
\end{enumerate}
\end{definition}The following examples illustrate
the notion of disjoint \emph{pieces}.

\begin{figure}[ht]
\centering
    \scalebox{0.4}{\includegraphics{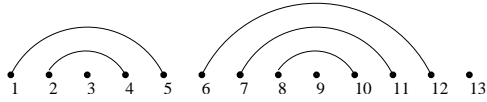}}
    \caption{$\pi$ made up of 3 disjoint pieces}
\end{figure}
\begin{figure}[ht]
\centering
    \scalebox{0.4}{\includegraphics{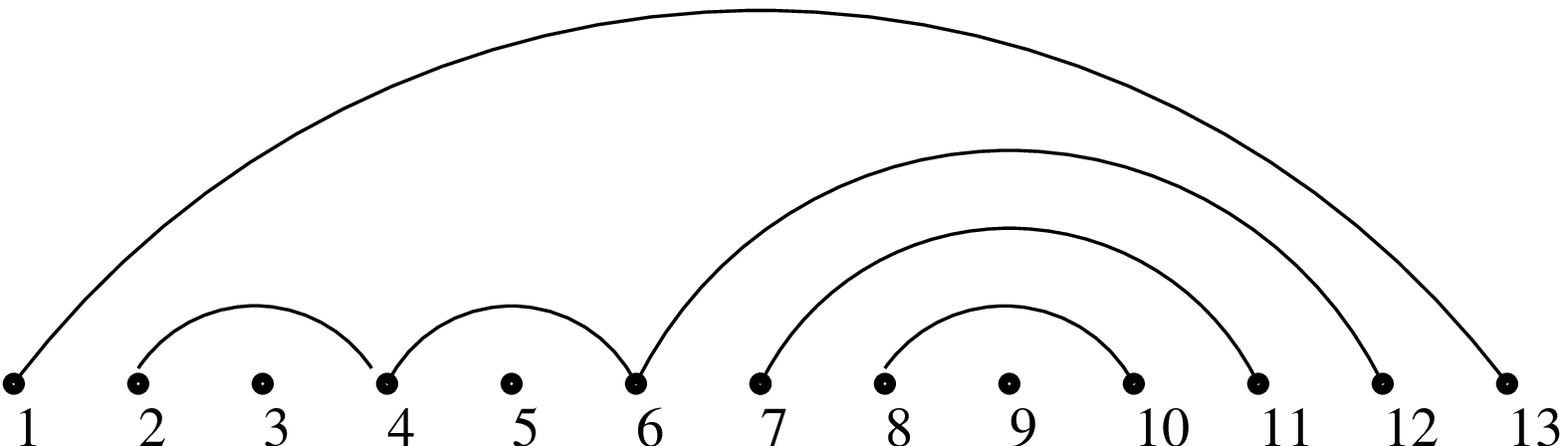}}
    \caption{$\pi$ made up of 1 piece}
\end{figure}

We also need the following lemma that will be used in the
proof of Theorem ~\ref{T1:theorem1}.

\begin{lemma}\label{L1:lemma1}
Let
$n,x_{1},x_{2},\dots,x_{k}$
be positive integers such that
$n=x_{1}+x_{2}+\dots+x_{k}$.
Then $\lfloor  \frac{x_{1}}{2} \rfloor +\lfloor  \frac{x_{2}}{2} \rfloor +
\dots+ \lfloor  \frac{x_{k}}{2} \rfloor + k \geq \lfloor  \frac{n}{2} \rfloor+1$,
where $\lfloor$ \ $\rfloor$ denotes the floor function.
\end{lemma}

\begin{proof}
First note that the case $k=1$ is trivially true.  So, let $k\geq2$ and
assume that $i$ of the $x_{j}$ are even, where $i \geq 0$.  Then,
\[
\sum_{j=1}^{k} \lfloor  \frac{x_{j}}{2} \rfloor+k
= \frac{n-(k-i)}{2} + k = \frac{n}{2} + \frac{k+i}{2} \geq \lfloor \frac{n}{2} \rfloor+1,
\forall k \geq 2.
\]
\end{proof}

\begin{theorem}\label{T1:theorem1}
Let $r$ be the number of blocks of a non-crossing partition of $[n]$
in which no block has two consecutive integers, (hereby referred to
as a semi-special partition, or, SSP).  Then
$r \geq \lfloor  \frac{n}{2} \rfloor + 1$.
\end{theorem}
\begin{proof}
We will use induction on $n$.
\begin{center}
\begin{enumerate}
    \item $n=1$.  The SSP of the singleton $[1]$
              contains exactly 1 block and clearly $1 \geq \lfloor
              \frac{1}{2} \rfloor + 1 = 1$.
    \item $n=2$.  The SSP  $\{\{1\},\{2\}\}$ of $[2]$
             contains two blocks, and again
             $2 \geq \lfloor \frac{2}{2} \rfloor +1=2$.

    \item   Assume the claim holds for all SSP's
        of $k \leq n$, and let $\pi$ be a SSP of $[n]$.  \\[.2cm]

\begin{enumerate}
    \item           Let $n+1=2j+1$.  By induction, $\pi$ is made up of
            $r \geq \lfloor  \frac{2j}{2} \rfloor + 1=j+1$ blocks.
            Let $\pi^{'}$ be a SSP of $[n+1]$, with $r^{'}$
            blocks.\\[.2cm]

            If $\pi^{'}$ is obtained from $\pi$ by adding
$n+1$ to an existing block, then  $r^{'} \geq j+1 = \lfloor  \frac{2j+1}{2} \rfloor + 1 =
            \lfloor  \frac{n+1}{2} \rfloor + 1$.\\[.2cm]

            Else, if $n+1$ is added as a singleton to $\pi$, then
            $r^{'}= r + 1 \geq j+2 \geq
            \lfloor  \frac{2j+1}{2} \rfloor + 1=\lfloor \frac{n+1}{2} \rfloor + 1$.
            Thus the claim holds for $n+1$ odd.\\[.2cm]

    \item           Let $n+1=2j+2$, and assume that $\pi^{'}$ is a SSP
of $[n+1]$ with $r^{'}$
            blocks.  We must show that $r^{'} \geq
            \lfloor  \frac{2j+2}{2} \rfloor + 1=j+2$.
            One sees that this inequality would fail to hold if
            $\pi^{'}$ is obtained
            by adding $n+1$ to an existing block of
            $\pi$
            with exactly $r=\lfloor  \frac{2j+1}{2} \rfloor +
1=j+1$ blocks.  In \cite{Mullin}, Mullin and Stanton proved that
$2n+1$ is the largest possible value of $k$ for which there exists
a SSP with $n+1$ blocks, thus showing that $n+1$ cannot be added
to
an existing block of $\pi$.  However, their result follows from a
theorem using face maps.  Here, we give a proof relying only on
elementary observations for SSP.  Moreover, this proof shows
that special partitions are made up of only one piece, a result which
is useful in the main theorem of section 3.
\\[.1in]

            We show that if $\pi$
            has exactly $j+1$ blocks, then the only SSP of $n+1$
            that one can obtain from it is by adding $n+1$ as a singleton.
            To prove this,
            we show that a SSP
            $\pi$ of $[n]=[2j+1]$ with exactly $j+1$ parts has
            exactly one piece in the sense of
            Definition 2.1.\\[.1in]

Assume that $\pi$ is made up of
exactly
$j+1$
blocks, which are split into 2 disjoint pieces, $B$ and $B^{*}$
(By Lemma ~\ref{L1:lemma1}, it suffices to consider the case where we have
2
disjoint pieces, as shown in Figure 4). Note that $B$ and $B^{*}$
are SSPs of $\pi$.
\begin{figure}[ht]
\centering
    \scalebox{0.8}{\includegraphics{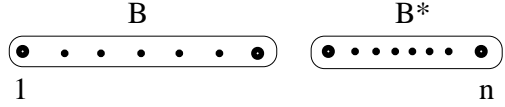}}
    \caption{$\pi$ made up of disjoint pieces}
\end{figure}
Since $n$ is odd we may assume that $B$
contains $2k$ integers, while $B^{*}$
contains $n-2k=2j+1-2k$ integers.
By induction, $B$ is made up of at least
$k+1$ blocks, while $B^{*}$
contains at least $j-k+1$ blocks.
Together, $B$ and $B^{*}$ contain
at least $j+2$ blocks, which contradicts our assumption on
$\pi$.
Thus $\pi$ is made up of 1 piece, and it follows that there
is exactly one way to add $n+1$ to $\pi$,  namely as a singleton.  Adding $n+1$ in any
other way will either violate the
non-crossing condition, or will cause the two consecutive integers $n$ and
$n+1$ to lie in the same block.
Therefore, $r^{'}= r + 1 \geq j+2 \geq \lfloor  \frac{2j+1+1}{2} \rfloor + 1 =\lfloor \frac{n+1}{2}\rfloor + 1$.

\end{enumerate}
\end{enumerate}
\end{center}

\end{proof}
The following corollary is an immediate consequence of
the above proof.
\begin{corollary}\label{C:corollary}
In a special partition of $[2n+1]$, $1$ and $2n+1$
are always in the same block.
\end{corollary}
\begin{proof}
Since a special partition consists of only 1 piece, $1$ and
$2n+1$ must be in the same block.
\end{proof}
\begin{corollary}\label{C1:corollary1}
The difference between any 2 consecutive
entries of any non-singleton block of a special partition
is always even.  That is, $b_{i_{j+1}}-b_{i_{j}}=2d$.
\end{corollary}

\begin{proof}
Consider a block $B_{i} \in \pi$, where $\pi$ is a special partition of $[2n+1]$.
Suppose there exists a difference of $2k+1$ between some integers
$b_{i_{j}}$ and $b_{i_{j+1}}$ in $B_{i}$.
It follows that there are $2k$ integers trapped between $b_{i_{j}}$ and
$b_{i_{j+1}}$, which
by Theorem ~\ref{T1:theorem1} make up at least $\lfloor
\frac{2k}{2} \rfloor+1=k+1$
blocks.
By identifying these $2k$ integers as
a $\textquoteleft$singleton',  we can view $\pi$ as a SSP of $[2n+1-2k+1]$
with at least
$\lfloor \frac{2n+1-2k+1}{2} \rfloor+1=n-k+2$ blocks.
Thus, $\pi$ is made up of at least $n-k+2+k+1-1=n+2$
blocks, a contradiction.
Therefore, the differences must be even.
\end{proof}

\begin{definition}
Given a block $B_{i}$ of size $m>1$ in a special partition $\pi$, a subpartition is a partition formed
by the set of integers between $b_{i_{j}}$ and $b_{i_{j+1}}$ of $B_{i}$, $1 \leq j < m$.
\end{definition}

\begin{corollary}\label{C2:corollary2}
Subpartitions are special partitions.
\end{corollary}
\begin{proof}
Given a special partition $\pi$ of $[2n+1]$, the subpartition
on $b_{i_{j}}+1, b_{i_{j}}+2, \dots, b_{i_{j+1}}-1$ is still non-crossing.
Moreover, any block within the subpartition does not contain two
consecutive integers.  It remains to show that for a difference of
$2k$ between $b_{i_{j}}$ and $b_{i_{j+1}}$, the $2k-1$ integers
$b_{i_{j}}+1, b_{i_{j}}+2, \dots, b_{i_{j+1}}-1$ are split into $k$
blocks.
This can be shown using the same technique as in the proof
of Corollary ~\ref{C1:corollary1}.
\end{proof}

\section{A Bijection Between Special Non-crossing Partitions and Sequences}
In this section we prove our main result, that is, we exhibit a bijection between the sets
$\mathcal P_{[2n+1]}$ and $ \mathcal S_{n}$, where $\mathcal S_{n}$ is the set of
sequences $s_{1}, s_{2}, \dots, s_{n}$ such that the following 2
conditions hold:
\begin{list}{}{\setlength{\leftmargin}{.5in}
        \setlength{\rightmargin}{.5in}}
    \item{i)\;} $1 \leq s_{i} \leq i$
    \item{ii)\;} if $s_{i}=j$, then $s_{i-r} \leq j-r$ for $1 \leq r \leq j-1$ \end{list} An element
$s \in \mathcal S_{6}$, for instance, is:  \begin{center} $s= 1 \; 2 \; 3 \; 1 \; 1 \; 6$ \end{center}
The following observations about sequences in $\mathcal S_{n}$ will
be useful: The
only integer
$s_{i}$ such that $s_{i}=s_{i+1}$
is $s_{i}=1$.

\begin{lemma}
Any sequence $s \in \mathcal S_{n}$ can be obtained
from a sequence of $n$ $1s$ using the following algorithm.
\begin{enumerate}
    \item Start with the two sequences $T$ and $g$, where
$T=1 \; 1 \; \dots \; 1$ (a sequence of
        $n$ 1's), and $g=1 \; 2 \; \dots \; n$.
    \item Build successive pairs of sequences
        $ \{T^{1},g^{1}\}, \{T^{2},g^{2}\},\dots,\{T^{n},g^{n}\}$ in the following
manner:

    \begin{enumerate}
            \item Pick an integer $m_{1} \leq n=g_{n}$.
            Let $T^{1}=T_{1}^{1} \; T_{2}^{1} \; \dots \;T_{n}^{1},$
            where $T_{i}^{1}=T_{i}$ for all $i$ such that $i<n$, and
            $T_{n}^{1}=m_{1}$.\\
            Let $g^{1}=g_{1}^{1} \; g_{2}^{1} \; \dots \; g_{n}^{1}$, where
            $g_{i}^{1}=g_{i}$ for all $i$ such that
            $i \leq n-m_{1}$, and $g_{n-m_{1}+j}^{1}=j$ for
            all $j$ such that $1 \leq j \leq m_{1}$.
            Note that choosing $m_{1}=1$ is also a possibility, in which case,
            $g^{1}=1 \; 2 \; \dots \; (n-1) \; 1$, and $T^{1}=T$.

            \item Build $T^{2}$ from $T^{1}$ by replacing the $(n-1)^{th}$ term
            of $T^{1}$, $T_{n-1}^{1}$, with any integer $m_{2} \leq g_{n-1}^{1}$.\\
            Let $g^{2}=1 \; 2 \; \dots \; n-m_{1} \; 1 \; 2 \; \dots \; (m_{1}-m_{2}-1) \;
            1 \; 2 \; \dots m_{2} \; m_{1}$. Note that
            if $m_{2} < m_{1}$, then
            $m_{1}-m_{2}-1 \geq 0$. In the case
            $m_{2}=m_{1}+1$, we simply have that $g^{2}=g^{1}$.

            \item Continue building such pairs of sequences $\{T^{i},g^{i}\}$
                where for all $i$, $T^{i}$ is obtained from $T^{i-1}$ by
                replacing the one in position $(n-i+1)$ with any integer
                $m_{i} \leq g_{n-i+1}^{i-1}$.  The sequence $g^{i}$ is
                built from the sequence $g^{i-1}$ in the manner described
                earlier. Notice that for all $i$, $T^{i},g^{i} \in \mathcal S_{n}$.
        \end{enumerate}
One clearly sees
that any sequence $s \in \mathcal S_{n}$ can be obtained in the manner described above,
and that indeed the only sequences obtained from this algorithm belong to $\mathcal S_{n}$.
\end{enumerate}
\end{lemma}

The following example illustrates this lemma.
\begin{example} Pairs of sequences
obtained from the algorithm above, where $n=8$.
    \begin{center}
        $T: 1 \; 1 \; 1 \; 1 \; 1 \; 1 \; 1 \; 1 \; \; \; g: 1
\; 2 \; 3 \; 4 \; 5 \; 6 \; 7 \; 8$\\

        $T^{1}: 1 \; 1 \; 1 \; 1 \; 1 \; 1 \; 1 \; 4 \; \; \;
g^{1}: 1 \; 2 \; 3 \; 4 \; 1 \; 2 \; 3 \; 4$\\

        $T^{2}: 1 \; 1 \; 1 \; 1 \; 1 \; 1 \; 1 \; 4 \; \; \;
g^{2}: 1 \; 2 \; 3 \; 4 \; 1 \; 2 \; 1 \; 4$\\

        $T^{3}: 1 \; 1 \; 1 \; 1 \; 1 \; 2 \; 1 \; 4 \; \; \;
g^{3}: 1 \; 2 \; 3 \; 4 \; 1 \; 2 \; 1 \; 4$\\

        $T^{4}: 1 \; 1 \; 1 \; 1 \; 1 \; 2 \; 1 \; 4 \; \; \;
g^{4}: 1 \; 2 \; 3 \; 4 \; 1 \; 2 \; 1 \; 4$\\

        $T^{5}: 1 \; 1 \; 1 \; 4 \; 1 \; 2 \; 1 \; 4 \; \; \;
g^{5}: 1 \; 2 \; 3 \; 4 \; 1 \; 2 \; 1 \; 4$\\[1.2cm]

    \end{center}
\end{example}

    Given $\pi \in \mathcal P_{[2n+1]}$, where $B_{i} = \{b_{i_{1}}, b_{i_{2}}, \dots, b_{i_{m}}\}$,
    we construct a sequence $d_{1},d_{2},\dots,d_{2n+1}$ of differences
    in the following manner:
\begin{enumerate}

\item Let $i=b_{j_{k}} \in B_{j}$.  Then let
\begin{list}{}{\setlength{\leftmargin}{.5in}
        \setlength{\rightmargin}{.5in}}
    \item{a)\;} $d_{i}=b_{j_{k+1}}-i$  if $|B_{j}| > k$.
    \item{b)\;} $d_{i}=0$ if $|B_{j}| = k$.
\end{list}
Since a special partition $\pi \in \mathcal P_{[2n+1]}$ is made up of $n+1$ blocks,
    there are exactly $n+1$ differences that are zero.  Thus, there
    are $2n+1-(n+1) = n$ non-zero differences.\\[.1cm]

    \item Let $\{d_{i_{k}}\}_{1}^{n}$ be the subsequence of non-zero
    integers,
    and form a new sequence $\{a_{j}\}_{1}^{n}$ such that
    $a_{j}=(d_{i_{n-j+1}})/2$.\\[.1cm]
\end{enumerate}

Define a map $f: \mathcal P_{[2n+1]} \mapsto [n-1]^{n}$ by
$f(\pi)=(a_{1}, a_{2}, \dots, a_{n})$,
where the sequence $(a_{1}, a_{2}, \dots, a_{n})$ is obtained using steps 1 and 2
as described above. \\[.3cm]
\begin{example}
\begin{center}
$\pi = \{1,13\} - \{2,4,6,12\} - \{3\} - \{5\} - \{7,11\} - \{8,10\} -  \{9\}$\\[.5cm]
\end{center}
The sequence of differences as obtained from step 1 is:
\begin{center}
 $d_{1}=12, d_{2}=2, d_{3}=0, d_{4}=2, d_{5}=0, d_{6}=6,$
\end{center}
\begin{center}
    $d_{7}=4, d_{8}=2, d_{9}=0, d_{10}=0, d_{11}=0, d_{12}=0, d_{13}=0$\\[.5cm]
\end{center}

The subsequence of non-zero integers is:
\begin{center}
$\{d_{i_{k}}\}_{1}^{6}=$  : $\{12 \quad 2 \quad 2 \quad 6 \quad 4 \quad 2$\} \\[.5cm]
\end{center}

The sequence after reversing and dividing by two is:
\begin{center}
$\{a_{j}\}_{1}^{6}=$: $\{1 \quad 2 \quad 3 \quad 1 \quad 1
 \quad 6\}$\\[.4cm]
\end{center}
\end{example}

It is easy to check that this sequence is an element of $\mathcal{S}_{6}$.
In Corollary ~\ref{C1:corollary1}, we showed that all differences $d_{i}$
are even.  Thus,
the sequence $\{a_{j}\}$ obtained from $f$ is indeed a sequence
of integers.

\begin{theorem}
The mapping $f$ is a bijection between $\mathcal{P}_{[2n+1]}$ and
$\mathcal{S}_{n}$.
\end{theorem}

\begin{proof}
First, we show that $f$ is well-defined.
The maximum value of $d_{1}$ is $2n+1-1=2n$.  But $(d_{1})/2$ corresponds
to $a_{n}$, so $a_{n} \leq n$.  Similarly,
the maximum value of $d_{i_{k}}$ is $2n-2k+2$.  Thus,
\[
a_{j}=d_{i_{n-j+1}}/2 \leq (2n-2(n-j+1)+2)/2 = 2j/2 = j,
\]
showing that the first condition for being a sequence in the set
$\mathcal {S}_{n}$ holds. By Corollary ~\ref{C2:corollary2}, any
difference of $2j$ between two integers $b_{i_{k}}$ and
$b_{i_{k+1}}$ of a block gives rise to a subpartition with $j-1$
differences. The maximum values of these differences are $2j-2,
2j-4, \dots, 2j-2(j-1)$. Therefore, $a_{j-r} \leq j-r$, $1 \leq r
\leq j-1$, and condition 2 holds. Thus, $f$ is well-defined.

To show that $f$ is onto, given a sequence $s=s_{1} s_{2} \dots
s_{n}
\in \mathcal S_{n}$, we build its corresponding special
partition in the following manner:\\[1.2cm]

Associate the sequence of $n$ $1s$
with the Puttenham diagram $\mathcal D_{1}$ below:

\begin{figure}[ht] \centering
    \scalebox{0.5}{\includegraphics{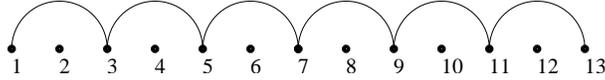}}
    \caption{$\mathcal D_{1}$: Puttenham diagram for a sequence with
$n$ $1's$, $n=6$}
\end{figure}

Note that $\mathcal D_{1}$
contains $n+1$ blocks and is clearly special.
We will construct a sequence of
special partitions from $\mathcal D_{1}$ as follows.
First, in order to obtain the sequence $s$ from the algorithm described in
Lemma 3.1, observe that $T_{n}^{1}=s_{n}$.  Now, create a corresponding
Puttenham diagram $\mathcal D_{2}$ by replacing the 1st arc of $\mathcal D_{1}$, $(1,3)$, with
the arc stretching from $1$ to $1+2s_{n}$, and by shifting $s_{n}-1$ arcs of length
two underneath it. In other words, the arcs $(3,5),(5,7),(7,9),\dots$, are replaced
by $(2,4),(4,6),(6,8)$ and so on.  Note that if $s_{n}=1$, then $\mathcal D_{2} =
\mathcal D_{1}$.  Observe that $\mathcal D_{2}$ contains
$n+1$ blocks, and is clearly special. Moreover, $T^{1}=1 \; 1 \; \dots \; s_{n}$
and $g^{1}=1 \; 1 \; \dots \; (n-s_{n}) \; 1 \; 2 \dots n$, which means that
for the next $s_{n}-1$ steps, we will be transforming $\mathcal D_{2}$ by working
\textquoteleft underneath' the arc $(1,1+2s_{n})$.

Next, $T^{2}=1 \; 1 \; \dots \; s_{n-1} \; s_{n}$, and we obtain $\mathcal D_{3}$ from
$\mathcal D_{2}$ by replacing the second arc of $\mathcal D_{2}$, $(2,4)$, with
the arc $(2,2+2s_{n-1})$ and once again shifting arcs of length two underneath
as before.  Again, if $s_{n-1}=1$, then $\mathcal D_{3}= \mathcal D_{2}$.

Continue creating each $\mathcal D_{i+1}$ in the following manner:
Stretch the $i$th arc of $\mathcal D_{i}$ by
$2(T_{n-i+1}^{i})=2s_{n-i+1}$, and
shift $T_{n-i+1}^{i}-1=s_{n-i+1}-1$ arcs \textquoteleft underneath' it.
Notice that after $s_{n}-1$ steps, we
repeat the process of stretching and shifting arcs at the integer $1+2s_{n}$, and so on.

The reason that we are able to shift arcs in the manner described above
is because each $T^{j}$ has a corresponding governing sequence $g^{j}$ whose
$(n-j+1)^{th}$ entry helps us identify the
exact number of arcs of length two that will
fit "underneath" the stretched arc in $\mathcal D_{j+1}$.
Also, notice that at each step, $\mathcal D_{j+1}$ is made up
of exactly one piece.

Clearly, each time arcs are stretched and shifted, the number of blocks in
the diagram $\mathcal D_{j+1}$ is
preserved, as are the conditions of being a special partition.
Note that subpartitions of each $\mathcal D_{j+1}$ are special
partitions as well. After $n$ steps (or less), this method of
construction gives us an element in $\mathcal P_{[2n+1]}$ that corresponds to $s$ and we are done.
The figure below illustrates how to construct the special partition corresponding to
$s= 1 \; 1 \; 1 \; 4 \; 1 \; 2 \; 1 \; 4$, where $n=8$. Notice here that
$T^{5}=T^{6}=T^{7}=T^{8}$,
and therefore $\mathcal D_{6} = \mathcal D_{7} = \mathcal D_{8} = \mathcal
D_{9}$.\\[5.0cm]

\begin{figure}[ht]
\centering
    \scalebox{0.5}{\includegraphics{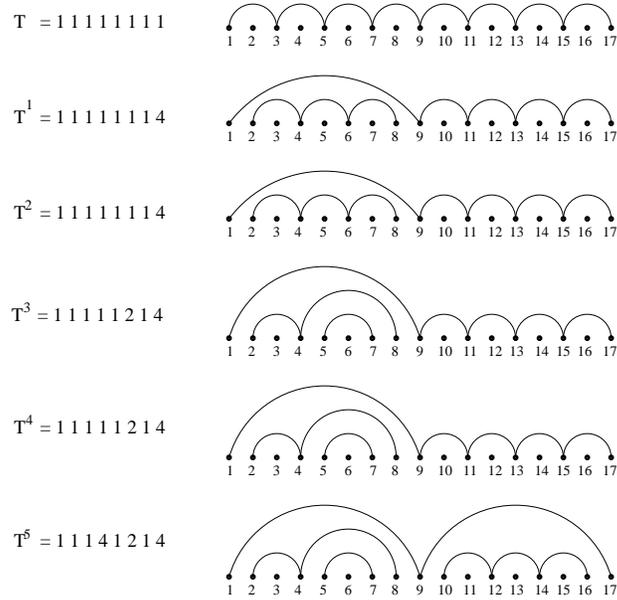}}
    \caption{Constructing the special partition corresponding to \; \;
        $s= 1 \; 1 \; 1 \; 4 \; 1 \; 2 \; 1 \; 4$.}
\end{figure}
\end{proof}

\section{Acknowledgments}
Sincere thanks are due to H\'{e}l\`{e}ne Barcelo for helpful
discussions and suggestions.

\end{document}